\documentclass[12pt]{article}

\usepackage{mathrsfs}

\usepackage{amssymb}
\usepackage{latexsym}
\usepackage{amsfonts}
\usepackage{amsmath}
\usepackage{eucal}
\usepackage{bm}
\usepackage{bbm}
\usepackage{graphicx}
\usepackage[english]{varioref}
\usepackage[nice]{nicefrac}
\usepackage[all]{xy}
\usepackage{amsthm}
\usepackage{amssymb,amsthm,upref,amscd}

\def\op{\operatorname}

\def\mmod{\kern-1pt\operatorname{-mod}}

\newtheorem{Thm}{Theorem}[section]
\newtheorem{Lem}[Thm]{Lemma}
\newtheorem{Cor}[Thm]{Corollary}
\newtheorem{Prop}[Thm]{Proposition}

\newtheorem{Rem}[Thm]{Remark}

\begin{document}
\renewcommand{\baselinestretch}{1.0}
\setlength{\parindent}{20pt} \baselineskip18pt

\centerline {\large \bf The Submodule Structure of the Permutation Module}
\centerline{\large \bf on Flag Varieties in Cross Characteristic}

\bigskip

\centerline{Xiaoyu Chen, Junbin Dong}

\begin{abstract}
Let ${\bf G}$ be a connected reductive group defined over $\mathbb{F}_q$ (the finite field with $q$ elements, where $q$ is a power of the prime number $p$), with the standard Frobenius map $F$. Let ${\bf B}$ be an $F$-stable Borel subgroup. Let $\Bbbk$ be a field (may not be algebraically closed) of characteristic $0\le r\neq p$. In this paper, we completely determine the composition factors of the permutation module module $\Bbbk[{\bf G}/{\bf B}]=\Bbbk{\bf G}\otimes_{\Bbbk{\bf B}}\op{tr}$ (here $\Bbbk{\bf H}$ is the group algebra of the group ${\bf H}$, and $\op{tr}$ is the trivial ${\bf B}$-module). In particular, we find a large family of infinite dimensional absolutely irreducible abstract representations of ${\bf G}$ over $\Bbbk$.
\end{abstract}

\section{Introduction}
The study of the induced modules of a reductive group ${\bf G}$ from a 1-dimensional representation of a subgroup is critical in understanding the irreducible modules of ${\bf G}$. In general, it is difficult to give such a decomposition. The well known Kazhdan-Lusztig (see. \cite{KL}) theory and Deligne-Lusztig theory (see. \cite{DL}) are closely related to certain type of induced modules: Verma modules and the Deligne-Lusztig's induction, respectively.

For the finite groups of Lie type, the submodule structure of such modules has been well studied (see.~\cite{C} and \cite{L}). For example, in \cite{J} Jantzen constructs a filtration for such induced modules and gives the sum formulas of these filtrations correspond to those of the well known Jantzen filtrations of generic Weyl modules. In \cite{P} C. Pillen proved that the socle and radical filtrations of such modules can be obtained from the filtrations of the generic Weyl modules under the similar assumption in \cite{J}. It was also showed in the same paper that these modules are rigid.

For the abstract representations of connected reductive groups, little was known for such modules. Recently, as observed by Nanhua Xi in \cite{X}, the induction using the construction of Frobenius, will produce new infinite dimensional irreducible abstract representations. For example, the infinite dimensional Steinberg module, which is the direct limit of Steinberg modules of $G_{q^a}$ (the $\mathbb{F}_{q^a}$ points of ${\bf G}$) as $a\rightarrow\infty$, was turned out to be irreducible in \cite{X} (in defining characteristic) and \cite{Y} (in cross characteristic). More generally, Xi constructed a filtration (whose subquotients are denoted by $E_J$, indexed by the subsets $J$ of the set $I$ of simple reflections) of the naive induced module from the trivial module of a Borel subgroup (see Section 2). It was proved in \cite{D} that some subquotients of this filtration are irreducible in cross characteristic.

The present paper shows that the above filtration given by Xi is a composition series. Thus, unlike the case of finite groups of Lie type, it is surprising that the naive induced module from the trivial module of a Borel subgroup has exactly $2^f$ composition factors, where $f$ is the rank of ${\bf G}$.  The idea of the proof contains the following three steps: (1) We show that $E_J$ is isomorphic to a submodule of certain parabolic induction; (2) Using the self duality of the finite version of such parabolic induction, we show that each submodule of $E_J$ contains a nonzero $B_{q^a}$-fixed point in $E_J$ for some $a>0$, where $B_{q^a}$ is the $\mathbb{F}_{q^a}$-points of ${\bf B}$; (3) We show that such $B_{q^a}$-fixed point can be translated to a generator of $E_J$ by $\Bbbk{\bf G}$. The above three steps imply the irreducibility of $E_J$. Xi proved that all these $E_J$ are pairwise non-isomorphic (see.~\cite[Proposition 2.7]{X}). We conclude that all $E_J$ with $J\subset I$ are irreducible and pairwise non-isomorphic.

 This paper is organized as follows: In Section 2 we recall some notation and basic facts. The Section 3 gives the proof of the main theorem, and gives some consequences of main theorem. In the section 4, we show that $E_J$ is not an union of irreducible $G_{q^a}$-modules (although ${\bf G}$ is union of $G_{q^a}$). The Section 5 gives some open problems for further study.

\medskip
\noindent{\bf Acknowledgements.} Xiaoyu Chen is supported by National Natural Science Foundation of China (Grant No.~11501546). Junbin Dong is supported by National Natural Science Foundation of China (Grant No.~11671297).  The authors would like to thank
Professor Nanhua Xi for his helpful suggestions and comments in
writing this paper.  The first author thanks Professor Jianpan Wang and Naihong Hu for their advice and comments.

\section{Notation and Preliminaries}
Let $\Bbbk$ be a field. Let ${\bf G}$ be a connected reductive group defined over $\mathbb{F}_q$ with the standard Frobenius map $F$. Let ${\bf B}$ be an $F$-stable Borel subgroup, and ${\bf T}$ be an $F$-stable maximal torus contained in ${\bf B}$, and ${\bf U}=R_u({\bf B})$ the ($F$-stable) unipotent radical of ${\bf B}$. We denote $\Phi=\Phi({\bf G};{\bf T})$ the corresponding root system, and $\Phi^+$ (resp. $\Phi^-$) be the set of positive (resp. negative) roots determined by ${\bf B}$. Let $W=N_{\bf G}({\bf T})/{\bf T}$ be the corresponding Weyl group. For each $w\in W$, let $\dot{w}$ be a representative in $N_{\bf G}({\bf T})$. For any $w\in W$, let ${\bf U}_w$ (resp. ${\bf U}_w'$) be the subgroup of ${\bf U}$ generated by all ${\bf U}_\alpha$ (the root subgroup of $\alpha\in\Phi^+$) with $w\alpha\in\Phi^-$ (resp. $w\alpha\in\Phi^+$). The multiplication map ${\bf U}_w\times{\bf U}_w'\rightarrow{\bf U}$ is a bijection (see \cite[Proposition 2.5.12]{C}). One denotes $\Delta=\{\alpha_i\mid i\in I\}$ the set of simple roots and $S=\{s_i\mid i\in I\}$ the corresponding simple reflections in $W$. For any $J\subset I$, let $W_J$ and ${\bf P}_J$ be the corresponding standard parabolic subgroup of $W$ and ${\bf G}$, respectively.

Let $\mathbb{M}(\op{tr})=\Bbbk{\bf G}\otimes_{\Bbbk{\bf B}}\op{tr}=\Bbbk[{\bf G}/{\bf B}]$, where $\op{tr}$ is the trivial ${\bf B}$-module, and call it the {\it permutation module on flag variety}. Let ${\bf 1}_{\op{tr}}$ be a nonzero element in $\op{tr}$. For convenience, we abbreviate $x\otimes{\bf 1}_{\op{tr}}\in\mathbb{M}(\op{tr})$ to $x{\bf 1}_{\op{tr}}$. For any $J\subset I$, set
$$\eta_J=\sum_{w\in W_J}(-1)^{\ell(w)}\dot{w}{\bf 1}_{\op{tr}},$$
and let $\mathbb{M}(\op{tr})_J=\Bbbk{\bf G}\eta_J$. It was proved in \cite{X} that
\begin{equation}\label{MJ=KUW}
\mathbb{M}(\op{tr})_J=\Bbbk {\bf U}W\eta_J
\end{equation}
Following \cite[2.6]{X}, we define
$$E_J=\mathbb{M}(\op{tr})_J/\mathbb{M}(\op{tr})_J',$$
where $\mathbb{M}(\op{tr})_J'$ is the sum of all $\mathbb{M}(\op{tr})_K$ with $J\subsetneq K$. The following lemma was proved in \cite{X}.

\begin{Lem}[{\cite[Proposition 2.7]{X}}]\label{EJ}
If $J$ and $K$ are different subsets of $S$, then $E_J$ and $E_K$ are not isomorphic.
\end{Lem}

 For each $w\in W$, let
$$C_w=\sum_{y\leq w}(-1)^{\ell(w)-\ell(y)}P_{y,w}(1)y\in\Bbbk W,$$
where $P_{y,w}$ are Kazhdan-Lusztig polynomials. According to \cite{KL}, the elements $C_w$ with $w\in W$ form a basis of $\Bbbk W$ and $C_{w_J}=(-1)^{\ell(w_J)}\eta_J$ for any $J\subset I$ by \cite[Lemma 2.6 (vi)]{KL}, where $w_J$ is the longest element in $W_J$.
For $w\in W$, set $\mathscr{R}(w)=\{i\in I\mid ws_i<w\}$. For any $J\subset I$, define
$$
\aligned
X_J &\ =\{x\in W\mid x~\op{has~minimal~length~in}~xW_J\};\\
Y_J &\ =\{x\in X_J\mid \mathscr{R}(xw_J)=J\}.
\endaligned
$$
The following lemma is well known and the proof can be found in \cite[Proposition 2.3]{X} (see also \cite{D}).
\begin{Lem}\label{xsty}
Let $u\in{\bf U}_{\alpha_i}\backslash\{1\}$ and $w\in X_J$. Then

\noindent $\op{(1)}$ There exists $x,y\in{\bf U}_{\alpha_i}\backslash\{1\}$ and $t\in{\bf T}$ such that $\dot{s_i}u\dot{s_i}^{-1}=x\dot{s_i}ty$;

\noindent $\op{(2)}$ If $ww_J< s_iww_J$, then $\dot{s_i}u\dot{w}\eta_J=\dot{s_i}\dot{w}\eta_J$;

\noindent $\op{(3)}$ If $s_iw< w$, then $\dot{s_i}u\dot{w}\eta_J=x\dot{w}\eta_J$, where $x$ is defined in $\op{(1)}$.

\noindent $\op{(4)}$ If $s_iw> w$ and $s_iww_J< ww_J$, then $\dot{s_i}u\dot{w}\eta_J=(x-1)\dot{w}\eta_J$, where $x$ is defined in $\op{(1)}$.
\end{Lem}
\noindent The following results was proved in \cite{D}.

\begin{Lem}[{\cite[Lemma 2.6]{D}}]\label{wCJ}
Let $J\subset I$. Then the following sets form a basis of $\Bbbk WC_{w_J}$:

\smallskip
\noindent $(1)$ $\{wC_{w_J}\mid w\in X_J\}$;

\smallskip
\noindent $(2)$ $\{C_{xw_J}\mid x\in X_J\}$;

\smallskip
\noindent $(3)$ $\{wC_{w_J}\mid w\in Y_J\}\cup\{C_{xw_J}\mid x\in X_J\backslash Y_J\}$.
\end{Lem}

\begin{Lem}[{\cite[Lemma 2.7]{D}}]\label{YJCJ}
For any $J\subset I$, denote by $C_J$ the image of $\eta_J$ in $E_J$. Then $E_J$ is spanned by all $u\dot{w}C_J$ with $w\in Y_J$ and $u\in{\bf U}_{w_Jw^{-1}}$.
\end{Lem}

\noindent For any finite subgroup $H$ of ${\bf G}$, let $\underline{H}:=\sum_{h\in H}h\in\Bbbk{G}$. For each $F$-stable subgroup ${\bf H}$ of ${\bf G}$, denote $H_{q^a}$ the $\mathbb{F}_{q^a}$-points of ${\bf H}$ (for example, $U_{q^a}$, $B_{q^a}$, $T_{q^a}$, $U_{w,q^a}$ etc.), then ${\bf H}$ is identified with its $\bar{\mathbb{F}}_q$-points and ${\bf H}=\displaystyle\bigcup_{a=1}^{\infty}H_{q^a}$.

\begin{Lem}[{\cite[Lemma 2.9]{D}}]\label{Yang}
Assume that $a>0$ and $\op{char}\Bbbk\neq\op{char}\bar{\mathbb{F}}_q$. Let $M$ be a $\Bbbk{\bf G}$-module and $0\neq\eta\in M^{\bf T}$. Then $\eta$ is contained in the submodule of $M$ generated by $\underline{U_{q^a}}\eta$.
\end{Lem}

If ${\bf G}$ and ${\bf T}$ in Lemma \ref{Yang} is replaced by $G_{q^a}$ and $T_{q^a}$, respectively, then the result no longer holds. For example, assume that $\op{char}\bar{\mathbb{F}}_q=p$, $\op{char}\Bbbk=r$, and $p\equiv1~(\op{mod}r)$. Let $M=\Bbbk G_{q^a}{\bf 1}_{\op{tr}}$ and $\eta=\sum_{w\in W}\dot{w}{\bf 1}_{\op{tr}}\in M^{T_{q^a}}$. Then
$$\underline{U_{q^a}}\eta=\sum_{w\in W}q^{a\ell(w_0w)}\underline{U_{w^{-1},q^a}}\dot{w}{\bf 1}_{\op{tr}}=\sum_{w\in W}\underline{U_{w^{-1},q^a}}\dot{w}{\bf 1}_{\op{tr}},$$
where $w_0$ is the longest element of $W$. It follows that $\underline{U_{q^a}}\eta$ generates a trivial $G_{q^a}$-submodule of $\Bbbk G_{q^a}\eta$.
 Therefore, Lemma \ref{Yang} displays a special feature of infinite dimensional representations of reductive
groups.

\section{The proof of main theorem}
The main theorem of this paper is the following
\begin{Thm}\label{main}
Assume that $0\le\op{char}\Bbbk\neq\op{char}\bar{\mathbb{F}}_q$. Then all modules $E_J$ $(J\subset I)$ are irreducible and pairwise non-isomorphic. In particular, $\mathbb{M}(\op{tr})$ has exactly $2^f$ composition factors, where $f$ is the rank of ${\bf G}$.
\end{Thm}

To prove Theorem \ref{main}, we start with the certain submodule of the  parabolic induction $\mathbb{M}_K=\op{Ind}_{{\bf P}_K}^{\bf G}\op{tr}_K=\Bbbk{\bf G}\otimes_{\Bbbk{\bf P}_K}\op{tr}_K$ for $K\subset I$, where $\op{tr}_K$ is the trivial ${\bf P}_K$-module.  Let ${\bf 1}_K$ be a nonzero element in $\op{tr}_K$. For convenience, we abbreviate $x\otimes{\bf 1}_K\in \mathbb{M}_K$ to $x{\bf 1}_K$ as before.

\begin{Prop}\label{basis}
Let $J\subset I$ and $J'=I\backslash J$. Set $E_J'$ be the submodule of $\mathbb{M}_{J'}$ generated by $D_J:=\sum_{w\in W_J}(-1)^{\ell(w)}\dot{w}{\bf 1}_{J'}$. Then $\{u\dot{w}D_J\mid w\in Y_J, u\in{\bf U}_{w_Jw^{-1}}\}$ forms a basis of $E_J'$. In particular, $E_J$ is a $($nonzero$)$ quotient of $E_J'$.
\end{Prop}
\begin{proof}
The proof is almost same as that in \cite[Lemma 2.7]{D}. Clearly, all $u\dot{w}D_J$ with $w\in Y_J$ and $u\in{\bf U}_{w_Jw^{-1}}$ are linearly independent. It remains to show that $E_J'$ is spanned by $\{u\dot{w}D_J\mid w\in Y_J, u\in{\bf U}_{w_Jw^{-1}}\}$.

For each $w\in W$, set $c_w=(-1)^{\ell(w)}C_w{\bf 1}_{J'}\in E_J'$. By \cite[Lemma 2.6 (vi)]{KL}, we have $c_{w_J}=D_J$. Since $\mathbb{M}_{J'}$ is a quotient of $\mathbb{M}(\op{tr})$, we have
\begin{equation}\label{kUWD}
E_J'=\Bbbk{\bf U}WD_J=\sum_{w\in Y_J}\Bbbk{\bf U}\dot{w}D_J+\sum_{x\in X_J\backslash Y_J}\Bbbk{\bf U}c_{xw_J}
\end{equation}
by (\ref{MJ=KUW}) and Lemma \ref{wCJ} (3). By (\ref{kUWD}), it remains to show that $c_{xw_J}=0$ for $x\in X_J\backslash Y_J$. In fact, let $J\subsetneq\mathscr{R}(xw_J)=K$. Then $xw_J=x'w_K$ for some $x'\in X_K$, and hence
$$c_{xw_J}=c_{x'w_K}=(-1)^{\ell(x'w_K)}C_{x'w_K}{\bf 1}_{J'}=h\sum_{w\in W_K}(-1)^{\ell(w)}\dot{w}{\bf 1}_{J'}=0$$
for some $h\in \Bbbk W$ by (1) and (2) of Lemma \ref{wCJ}.

It follows that there is a surjective ${\bf G}$-module homomorphism $E_J'\rightarrow E_J$ mapping $u\dot{w}D_J$ to $u\dot{w}C_J$ by Lemma \ref{YJCJ}. This completes the proof.
\end{proof}

\bigskip
Before proving Theorem \ref{main}, we make some preliminaries. For each $i\in I$, let $u\in{\bf U}_{\alpha_i}\backslash\{1\}$, and $x\in{\bf U}_{\alpha_i}\backslash\{1\}$ be the elements defined in (1) of Lemma \ref{xsty}. Set $\tau_i:=u^{-1}\dot{s_i}^{-1}(x-1)\in\Bbbk{\bf G}$. We have the following trivial observation.
\begin{Lem}\label{tau}
For any $i\in I$ and $w\in Y_J$, we have
$$
\tau_i\dot{w}D_J=\left\{
\begin{array}{ll}
 \dot{w}D_J-\dot{s_i}\dot{w}D_J &\ \mbox{if}~s_iw<w\\
\dot{w}D_J &\ \mbox{if}~s_iww_J<ww_J~\mbox{and}~s_iw>w\\
 0 &\ \mbox{if}~s_iww_J>ww_J
\end{array}\right.
$$
\end{Lem}
\begin{proof}
The second case follows immediately from (4) of Lemma \ref{xsty}.
The third case is trivial since $x\dot{w}D_J=\dot{w}D_J$ if $s_iww_J>ww_J$.

If $s_iw<w$, then $u^{-1}\dot{s_i}\dot{w}D_J=\dot{s_i}\dot{w}D_J$. It follows from (3) of Lemma \ref{xsty} that
$$\tau_i\dot{w}D_J=u^{-1}\dot{s_i}^{-1}x\dot{w}D_J-\dot{s_i}\dot{w}D_J=\dot{w}D_J-\dot{s_i}\dot{w}D_J.$$
This completes the proof.
\end{proof}

\bigskip
\begin{Rem}
From Lemma \ref{tau}, we see that as an operator on $\Bbbk WD_J$, $\tau_i$ is independent of the choice of $u$ and $x$.
\end{Rem}

\begin{Cor}\label{tau2}
Let $j_1,\cdots,j_k\in I$. If the coefficient $($in terms of the basis given in Proposition \ref{basis}$)$ of $\dot{w_1}D_J$ in $\tau_{j_k}\cdots\tau_{j_1}\dot{w_2}D_J$ is nonzero, then $w_1=w_2$, or there exists a $1\leq t\leq k$ and a subset $\{i(1),i(2),\cdots,i(t)\}$ of $\{1,2,\cdots,k\}$ such that $(\op{i})$ $i(1)<i(2)<\cdots<i(t)$, $(\op{ii})$ $\ell(w_1)=\ell(w_2)-t$, and $(\op{iii})$ $w_1=s_{j_{i(t)}}\cdots s_{j_{i(1)}}w_2$.
\end{Cor}
\begin{proof}
We proceed by the induction on $k$. The case $k=1$ is trivial by Lemma \ref{tau}. Assume that $k>1$ and
$$\tau_{j_{k-1}}\cdots\tau_{j_1}\dot{w_2}D_J=\sum_{w\in Y_J}a_w\dot{w}D_J,\quad a_w\in\Bbbk.$$
Let $Y=\{w\in Y_J\mid a_w\neq0\}$ and assume that the coefficient of $\dot{w_1}D_J$ in $\tau_{j_k}\cdots\tau_{j_1}\dot{w_2}D_J$ is nonzero. Then the coefficient of $\dot{w_1}D_J$ in $\tau_{j_k}\dot{w'}D_J$ is nonzero for some $w'\in Y$. By induction we have $w'=w_2$, or there is some $1\leq t\leq k-1$ and a subset $\{i(1),i(2),\cdots,i(t)\}$ of $\{1,2,\cdots,k-1\}$ satisfying (i),(ii),(iii) with $w_1$ replaced by $w'$. By Lemma \ref{tau}, we have $w_1=w'$ or $w_1=s_{j_k}w'<w'$. Thus $w_1=w_2$, or we choose the subset $\{i(1),i(2),\cdots,i(t)\}$ in the former case, and $\{i(1),i(2),\cdots,i(t),k\}$ in the latter case. Clearly, these subsets satisfy (i), (ii), (iii) for $w_1$. This completes the proof.
\end{proof}

\bigskip
\noindent As an easy consequence of Corollary \ref{tau2}, we have
\begin{Cor}\label{tau3}
Let $j_1,\cdots,j_k\in I$. Then

\smallskip
\noindent $(1)$ The coefficient of $\dot{w_1}D_J$ in $\tau_{j_k}\cdots\tau_{j_1}\dot{w_2}D_J$ is zero if $\ell(w_2)-\ell(w_1)>k$.

\smallskip
\noindent $(2)$ If $\ell(w_2)-\ell(w_1)=k$, then the coefficient of $\dot{w_1}D_J$ in $\tau_{j_k}\cdots\tau_{j_1}\dot{w_2}D_J$ is nonzero if and only if $w_1=s_{j_k}\cdots s_{j_1}w_2$.
\end{Cor}
\begin{proof}
(1) and the ``only if" part of (2) is obvious by Corollary \ref{tau2}. An easy calculation using Lemma \ref{tau} shows that if $w_1=s_{j_k}\cdots s_{j_1}w_2$, then the above coefficient is $(-1)^k$, and hence the ``if" part of (2) holds. This completes the proof.
\end{proof}

\bigskip
Now we are going to prove the main result Theorem \ref{main}. From here to the end of this section, we assume that $\op{char}\Bbbk\neq\op{char}\bar{\mathbb{F}}_q$. We first prove the following

\begin{Prop}\label{submain}
For any $J\subset I$, the module $E_J'$ is irreducible.
\end{Prop}
\noindent To see this, we first prove the following two claims.
Throughout the proof, we assume that all representatives of the elements of $W$ involved are in $G_q$ without loss of generality. Otherwise we replace $q$ by a sufficiently large power of $q$. This does no harm to the proof.

\smallskip
\noindent {\bf Claim 1}. {\it If $M$ is a nonzero submodule of $E_J'$, then $M\cap\Bbbk WD_J\neq0$.}

\smallskip
\noindent {\bf Claim 2}. {\it If $N$ is a nonzero submodule of $E_J'$ such that $N\cap\Bbbk WD_J\neq0$, then $D_J\in N$ and hence $N=E_J'$.}

\smallskip
\noindent{\it Proof of claim 1.}
Assume that $M$ is a nonzero submodule of $E_J'$ and $0\neq x\in M$. Then $x\in E_{J,q^a}'=\Bbbk G_{q^a}D_J$ for some $a>0$.

\smallskip
For any $K\subset I$, set $\mathbb{M}_{K,q^a}=\op{Ind}_{P_{K,q^a}}^{G_{q^a}}\op{tr}_K$. Let $\mathcal{H}$ and $\mathcal{S}$ be the head and the socle of $\mathbb{M}_{J',q^a}$, respectively. By the self-duality of $\mathbb{M}_{J',q^a}$ (this follows from the well known result that the permutation module $P=\Bbbk[X]$ associated to a action of a finite group $G$ to a finite set $X$ is self dual, and the corresponding isomorphism sends each $x\in X$ to $f_x$, the characteristic function of $x$), $V\in\op{Irr}_\Bbbk(G_{q^a})$ embeds in $\mathcal{S}$ if and only if $V^*$ embeds in $\mathcal{H}$. Assume that $L\in\op{Irr}_\Bbbk(G_{q^a})$ embeds in $\mathcal{S}$. Then $$\op{Hom}_{P_{J',q^a}}(\op{tr}_{J'},L^*)=\op{Hom}_{G_{q^a}}(\mathbb{M}_{J',q^a},L^*)\neq0$$ by Frobenius reciprocity and the above discussion. This is equivalent to $(L^*)^{P_{J',q^a}}\neq0$, and hence $(L^*)^{B_{q^a}}\neq0$ which is equivalent to $L^{B_{q^a}}\neq0$ by {\cite[Lemma 2.1]{GP}}.

\smallskip
It is clear that $(E_{J,q^a}')^{B_{q^a}}\subset\bigoplus_{w\in Y_J}\Bbbk\underline{U_{w_Jw^{-1},q^a}}\dot{w}D_J$ by Proposition \ref{basis}, and there is an $L\in\op{Irr}_\Bbbk(G_{q^a})$ such that $L\subset\Bbbk G_{q^a}x\cap\mathcal{S}\subset E_{J,q^a}'\cap M$. By the previous paragraph, $L^{B_{q^a}}\neq0$, which implies that $(E_{J,q^a}')^{B_{q^a}}\cap M\neq0$. Assume that
$$0\neq\xi=\sum_{w\in Y_J}c_w\underline{U_{w_Jw^{-1},q^a}}\dot{w}D_J\in(E_{J,q^a}')^{B_{q^a}}\cap M,\quad c_w\in\Bbbk.$$
Notice that
$$
\aligned
\underline{U_{q^a}}\xi &\ =\underline{U_{q^a}}\sum_{w\in Y_J}c_w\underline{U_{w_Jw^{-1},q^a}}\dot{w}D_J\\
&\ =\sum_{w\in Y_J}\underline{U_{w_Jw^{-1},q^a}'}\cdot\underline{U_{w_Jw^{-1},q^a}}c_w\underline{U_{w_Jw^{-1},q^a}}\dot{w}D_J\\
&\ =\sum_{w\in Y_J}\underline{U_{w_Jw^{-1},q^a}'}\cdot\underline{U_{w_Jw^{-1},q^a}}c_wq^{a\ell(w_Jw^{-1})}\dot{w}D_J\\
&\ =\underline{U_{q^a}}\sum_{w\in Y_J}c_wq^{a\ell(w_Jw^{-1})}\dot{w}D_J\\
&\ \in M.
\endaligned
$$
Since $q\neq 0$ in $\Bbbk$, we see that $0\neq\sum_{w\in Y_J}c_wq^{a\ell(w_Jw^{-1})}\dot{w}D_J\in (E_J')^{\bf T}$. It follows that $\sum_{w\in Y_J}c_wq^{a\ell(w_Jw^{-1})}\dot{w}D_J\in M\cap\Bbbk WD_J$ by Lemma \ref{Yang}. This completes the proof.

\bigskip
\noindent{\it Proof of claim 2.}
Let $\preceq$ be the natural dictionary order on $\mathbb{Z}_{\geq0}\times\mathbb{Z}_{>0}$, that is,
$$(m,n)\preceq(m',n')\Leftrightarrow m<m'~\vee~(m=m'~\wedge~n\leq n').$$
For any $0\neq Z=\sum_{w\in Y_J}z_w\dot{w}D_J\in\Bbbk WD_J$, let
$$h(Z)=\op{max}(\{\ell(w)\mid w\in Y_J,z_w\neq0\})$$
and $c(Z)$ be the number of elements $w\in Y_J$ such that $\ell(w)=h(Z)$. Set $\Psi(Z)=(h(Z),c(Z))\in\mathbb{Z}_{\geq0}\times\mathbb{Z}_{>0}$.

\smallskip
Assume that $0\neq A=\sum_{w\in Y_J}a_w\dot{w}D_J\in N\cap\Bbbk WD_J$. We will prove Claim 2 by the induction on $\Psi(A)$ with respect to $\preceq$.

\smallskip
If $\Psi(A)=(0,1)$, then $A=kD_J$ for some $k\in\Bbbk^\times$, and hence $D_J\in N$. Now we assume that $\Psi(A)\succ(0,1)$ (so that $h(A)>0$). Choose a $w\in Y_J$ with $\ell(w)=h(A)$.  Let $w=s_{i_1}\cdots s_{i_t}$ be its reduced expression ($t=h(A)$). Let $K$ be the minimal subset of $I$ containing $J$ such that $ww_J<w_K$. Since $w\in Y_J$, $ww_J\neq w_K$, and hence $s_jww_J>ww_J$ for some $j\in K$. It follows that
\begin{equation}\label{eq1}
\tau_j\dot{w}D_J=0
\end{equation}
by Lemma \ref{tau}. By the minimality of $K$, we have
\begin{equation}\label{eq2}
s_js_{i_{t'}}\cdots s_{i_1}ww_J<s_{i_{t'}}\cdots s_{i_1}ww_J
\end{equation}
for some $1\leq t'\leq t$. We abbreviate $w_1=s_{i_{t'}}\cdots s_{i_1}w\in Y_J$ for convenience. Consider the following two cases.

\medskip
\noindent Case (i): $a_{w_1}\neq0$. It follows from (\ref{eq2}) that $(x-1)\dot{w_1}D_J\neq0$ if $x\in{\bf U}_{\alpha_j}\backslash\{1\}$, and hence $(x-1)A\neq0$ by Proposition \ref{basis}. Equivalently, we have $\tau_jA\neq0$. Combining Lemma \ref{tau} and (\ref{eq1}), we see that $\tau_jA\in N$ and $\Psi(\tau_jA)\prec\Psi(A)$. Applying the inductive hypothesis to $\tau_jA$ yields the result.

\medskip
\noindent Case (ii): $a_{w_1}=0$. Let $l\geq1$ be the minimal number with the following property:

\smallskip
\noindent {\it There exists $w_2\in Y_J$ such that $(\op{i})$ $a_{w_2}\neq0$; $(\op{ii})$ $\ell(w_2)=\ell(w_1)+l$; $(\op{iii})$ $w_1=s_{i_{n(l)}}\cdots s_{i_{n(1)}}w_2$ for some subset $\{n(1),n(2),\cdots n(l)\}$ of $\{1,2,\cdots,t'\}$ with $n(1)<n(2)<\cdots<n(l)$.}

\smallskip
\noindent Since $t'$ satisfies this property, we have $l\leq t'$. For $w'\in Y_J$, denote $\kappa(w_1,w')$ the coefficient of $\dot{w_1}D_J$ in $\tau_{i_{n(l)}}\cdots \tau_{i_{n(1)}}\dot{w'}D_J$. We have the following facts

\noindent $(1)$ If $\ell(w')>\ell(w_1)+l$, then  $\kappa(w_1,w')=0$;

\noindent $(2)$ If $\ell(w')=\ell(w_1)+l$, then  $\kappa(w_1,w')\neq0$ if and only if $w'=w_2$;

\noindent $(3)$ If $\ell(w')<\ell(w_1)+l$ and $a_{w'}\neq0$, then $\kappa(w_1,w')=0$;

\noindent Indeed, (1) and (2) follows immediately from Corollary \ref{tau3}, and (3) follows from Corollary \ref{tau2}, the minimality of $l$, and Lemma \ref{tau}. Combining the above three facts, we see that the coefficient of $\dot{w_1}D_J$ in $B=\tau_{i_{n(l)}}\cdots\tau_{i_{n(1)}}A\in N$ equals to that in $\tau_{i_{n(l)}}\cdots\tau_{i_{n(1)}}a_{w_2}\dot{w_2}D_J$ which is $(-1)^la_{w_2}\neq0$. It follows that $B\neq0$ and $\Psi(B)\preceq\Psi(A)$ by Lemma \ref{tau}. If $\Psi(B)\prec\Psi(A)$, then the result follows from the induction. If $\Psi(B)=\Psi(A)$, then the coefficient of $\dot{w}D_J$ in $B$ is nonzero. We choose $j\in K$ as before and the result follows from applying Case (i) to $B$. This completes the proof.

\bigskip
\noindent{\it Proof of Proposition \ref{submain} and Theorem \ref{main}.}
Combining Claim 1 and Claim 2 above, we see that any nonzero submodule $M$ of $E_J'$ contains $D_J$, and hence $M=E_J'$. In particular, all $E_J'$ are irreducible. Therefore, all $E_J$ are irreducible and pairwise nonisomorphic by Proposition \ref{EJ} and \ref{basis}. This completes the proof.

Since \cite[Lemma 2.1]{GP} holds without the assumption that $\Bbbk$ is algebraically closed (see its proof), the whole proof in this paper doesn't involve this assumption. Therefore, the irreducibility of $E_J$ holds for any field $\Bbbk$ with $\op{char}\Bbbk\neq\op{char}\bar{\mathbb{F}}_q$.

\begin{Cor}
The composition factors of $\mathbb{M}_J$ are exactly all $E_K$ with $K\cap J=\emptyset$.
\end{Cor}
\begin{proof}
Let $N_J=\sum_{i\in J}\mathbb{M}(\op{tr})_{\{s_i\}}$. Then
$\mathbb{M}_J=\mathbb{M}(\op{tr})/N_J$ by \cite[Theorem 6.3]{D}, and $\mathbb{M}(\op{tr})_K\nsubseteq N_J$ if and only if $K\cap J=\emptyset$. For each submodule $N$ of $\mathbb{M}(\op{tr})$, let $\overline{N}$ be its image in $\mathbb{M}_J$. There is a natural surjection $$E_K\rightarrow\overline{\mathbb{M}(\op{tr})_K}/\sum_{K\subsetneq L}\overline{\mathbb{M}(\op{tr})_L}$$
for each $K\subset I$ with $K\cap J=\emptyset$, and hence the result follows immediately from Theorem \ref{main}.
\end{proof}

\bigskip
\begin{Cor}
Let $V$ be an $($abstract$)$ irreducible representation of ${\bf G}$ and $V^{\bf B}\neq0$. Then $V$ is the trivial representation.
\end{Cor}
\begin{proof}
By the Frobenius reciprocity, $V^{\bf B}\neq0$ implies that $V$ is a quotient of $\mathbb{M}(\op{tr})$. By Theorem \ref{main}, we have $V=E_J$ for some $J\subset I$. It remains to show that $J=\emptyset$. Let $0\neq v\in (E_J)^{\bf B}$. Then $v\in\Bbbk G_{q^a}D_J$ for some $a>0$. It follows that $v\in(\Bbbk G_{q^a}D_J)^{B_{q^a}}$, and hence $v=\sum_{w\in Y_J}c_w\underline{U_{w_Jw^{-1},q^a}}\dot{w}D_J$ with $c_w\in\Bbbk$ and $c_{w'}\neq0$ for some $w'\in Y_J$. Suppose that $J\neq\emptyset$, then $w_Jw'^{-1}\neq1$, and hence ${\bf U}_{w_Jw'^{-1}}\neq\{1\}$. It follows that if $g\in {\bf U}_{w_Jw'^{-1}}\backslash U_{w_Jw'^{-1},q^a}$, then  $g\underline{U_{w_Jw'^{-1},q^a}}\dot{w'}D_J\neq\underline{U_{w_Jw'^{-1},q^q}}\dot{w'}D_J$, and hence $gv\neq v$ by Proposition \ref{basis} and Theorem \ref{main}, which contradicts to $v\in (E_J)^{\bf B}$. This completes the proof.
\end{proof}

\section{Non quasi-finiteness of $E_J$}

Firstly we recall the quasi-finiteness defined in \cite{X}. A group $H$ is quasi-finite if $H$ has a sequence $H_1$, $H_2$,$\cdots$, $H_n$, $\cdots$ of finite subgroups such that $H$ is the union of all $H_i$ and for any positive integers $i$, $j$ there exists integer $r$ such that $H_i$ and $H_j$ are contained in $H_r$. The sequence $H_1$, $H_2$,$\cdots$, $H_n$, $\cdots$ is called a quasi-finite sequence of $H$. An irreducible  representation $M$ of $H$ is called {\it quasi-finite} (with respect to the quasi-finite sequence $H_1$, $H_2$,$\cdots$, $H_n$, $\cdots$) if it has a sequence of subspaces $M_1$, $M_2$,$\cdots$, $M_n$, $\cdots$ such that (1) $M_i$ is an irreducible $H_i$-module, (2) if $H_i$ is a subgroup
of $H_j$, then $M_i$ is a subspace of $M_j$, and (3) $M$ is the union of all $M_i$. Clearly, the sequence $G_q$, $G_{q^2}$, $G_{q^3}$, $\cdots$ is a quasi-finite sequence of ${\bf G}$ and we fix such quasi-finite sequence in the following.

Since ${\bf G}=\displaystyle\bigcup_{a=1}^{\infty}G_{q^a}$, for each countable dimensional ${\bf G}$-module $M$ we have $M=\displaystyle\bigcup_{a=1}^{\infty}M_a$, where $M_a$ is a finite dimensional $G_{q^a}$-module, and $M_a\subset M_b$ if $a|b$. In particular, the above property holds for irreducible ${\bf G}$-modules. Xi proved that if each $M_a$ is irreducible, then $M$ is irreducible (\cite[Lemma 1.6 (b)]{X}). However, the converse is not true in general. The preprint \cite{Y} told  us that the Steinberg module $E_I$ is not quasi-finite irreducible when $0\ne\op{char}\Bbbk\neq\op{char}\bar{\mathbb{F}}_q$.

One can also define the ``finite version" of $E_J$ in the similar fashion. Specifically, let $M(\op{tr})_{J,q^a}=\Bbbk G_{q^a}\eta_J$. Let $E_{J,q^a}=M(\op{tr})_{J,q^a}/M(\op{tr})_{J,q^a}'$, where $M(\op{tr})_{J,q^a}'$ is the sum of all $M(\op{tr})_{K,q^a}$ with $J\subsetneq K$. In contrast to Theorem \ref{main}, $E_{J,q^a}$ is reducible in general. For example, $E_{I,q^a}$ is the ordinary Steinberg module of $G_{q^a}$, which may reducible in the cross characteristic.
The following proposition shows that $E_J$ is not quasi-finite irreducible (with respect to the sequence $G_q$, $G_{q^2}$, $G_{q^3}$, $\cdots$) in general.

\begin{Prop}\label{nonqf}
Assume that the length of the regular module $\mathbb{C}W$ is not equal to $2^{|I|}$.  Then $E_J$ is not quasi-finite irreducible (with respect to the sequence $G_q$, $G_{q^2}$, $G_{q^3}$, $\cdots$) for some $J\subset I$.
\end{Prop}

\begin{proof}
Suppose that all $E_J$ are quasi-finite. For each $J\subset I$, write $E_J=\displaystyle\bigcup_{a=1}^{\infty}V_a$, where each $V_a$ is an irreducible $\mathbb{C} G_{q^a}$-module and $V_m\subset V_n$ if $m|n$. Let $0\ne\xi\in V_1$. The irreducibility of each $V_a$ yields $V_a= \mathbb{C} G_{q^a}\xi$. Thus $E_J=\mathbb{C} {\bf G}\xi$ and $C_J\in \mathbb{C} G_{q^{a(J)}}\xi$ for some integer $a(J)$. The irreducibility of $V_{a(J)}$ implies $E_{J,q^{a(J)}}=V_{a(J)}$ and hence $E_{J,q^{a(J)}}$ is irreducible. Choose $c$ so that $a(J)|c$ for all $J\subset I$, and hence $E_{J,q^c}=V_c$ is irreducible for all $J\subset I$. Thus, the length of $\op{Ind}_{B_{q^c}}^{G_{q^c}}\op{tr}$ is $2^{|I|}$. However it is known that there is a bijection between the composition
factors of $\mathbb{C} G_{q^c}$-module $\op{Ind}_{B_{q^c}}^{G_{q^c}}\op{tr}$ and
the composition factors of the regular module $\mathbb{C} W$ of $W$, which
preserves multiplicities.  This contradicts to the assumption and the theorem is proved.
\end{proof}

\noindent{\bf Remark 4.2.} The examples satisfying the assumption in Proposition \ref{nonqf} are ubiquitous. We consider ${\bf G}=SL_n(\bar{\mathbb{F}}_q)$ and $W$ is the symmetric group $\mathfrak{S}_n$. Let $T(n)$ be the length of $\mathbb{C}\mathfrak{S}_n$, which is the number of standard Young tableaux with $n$ cells. Such number is also called {\it $n$-th telephone number} historically.
It is well known that the telephone numbers satisfy the recurrence relation
$$T(n)=T(n-1)+(n-1)T(n-2)$$
with $T(1)=1$ and $T(2)=2$. It follows easily from induction that $T(n)> 2^{n-1}=2^{|I|}$ if $n\ge 4$, which satisfies the assumption in Proposition \ref{nonqf}.

\section{Further questions}
In this section we propose some questions on infinite dimensional abstract representations of reductive groups with Frobenius maps. Any 1-dimensional representation $\theta$ of ${\bf T}$ is regarded as a representation of ${\bf B}$ through the homomorphism ${\bf B}\rightarrow{\bf T}$. Let $\mathbb{M}(\theta)=\Bbbk{\bf G}\otimes_{\Bbbk{\bf B}}\theta$. The following questions naturally arises.

\smallskip
\noindent $\bullet$ Assume that $\op{char}\Bbbk=\op{char}\bar{\mathbb{F}}_q$ and $\theta$ is a rational character of ${\bf T}$. What is the necessary and sufficient condition for $\mathbb{M}(\theta)$ to have finitely many composition factors?

\smallskip
\noindent $\bullet$ Assume that $\op{char}\Bbbk\neq\op{char}\bar{\mathbb{F}}_q$ and the stabilizer $W_\theta$ of $\theta$ in $W$ is a parabolic subgroup $W_J$ for some $J\subset I$. How to decompose $\mathbb{M}(\theta)$?

\smallskip
\noindent $\bullet$ Are all $E_J$ irreducible if $\op{char}\Bbbk=\op{char}\bar{\mathbb{F}}_q$?

\smallskip
\noindent $\bullet$ It is well known that the decomposition of $\Bbbk G_{q^a}{\bf 1}_{\op{tr}}$ is closely related to the Hecke algebras. How to develop a parallel theory for ${\bf G}$?

{\small
\noindent Xiaoyu Chen

\noindent E-mail: gauss\_1024@126.com

\noindent Affiliation: 1. Department of Mathematics, Ningbo University,
818 Fenghua Road, Ningbo 315211, PR China;

\noindent 2. Institute of Mathematics, Academy of Mathematics and Systems Science, Chinese Academy of Sciences, Beijing 100190, China.

\bigskip
\noindent Junbin Dong

\noindent E-mail: dongjunbin1990@126.com

\noindent  Affiliation: 1. Department of Mathematics, Tongji University
1239 Siping Road, Shanghai 200092, P. R. China;

\noindent 2. Institute of Mathematics, Academy of Mathematics and Systems Science, Chinese Academy of Sciences, Beijing 100190, China.
}
\end{document}